\documentclass[11pt]{article}
\usepackage{tikz}

\usepackage{amsmath,amssymb}
\usepackage{amsfonts}
\usepackage{epsfig}
\usepackage{graphics}
\usepackage{graphicx}
\usepackage[mathscr]{eucal}

\usepackage[all]{xy}
\usepackage{bm}

\usepackage{multicol}
\usepackage{multirow}
\usepackage{subfig}

\newtheorem{theorem}{Theorem}

\newtheorem{algorithm}[theorem]{Algorithm}
\newtheorem{lemma}[theorem]{Lemma}

\def\proof{\par{\bf Proof}. \ignorespaces}

\def\endproof{\vbox{\hrule height0.6pt\hbox{%
                    \vrule height1.3ex width0.6pt\hskip0.8ex
                    \vrule width0.6pt}\hrule height0.6pt
                   }}

\def\R {{\mathbb R}}

\def\cU {{\cal U}}

\title{\bf{Numerical solution for Fokker-Planck equation using a two-level scheme}}

\author{M.M. Butt\thanks{
Department of Mathematics and Statistics, King Fahd University of Petroleum and Minerals, Dhahran $31261$, Saudi Arabia.
({\tt mmunirbutt@gmail.com})}
}

\date{}
\begin{document}
\maketitle

\begin{abstract}
 A numerical solution to the Fokker-Planck equation using a two-level scheme is presented. The Fokker-Planck (FP) equation is of parabolic type equation govern the time evolution of probability density function of the stochastic processes. The FP equation also preserves the positivity and conservative of the total probability. A Chang-Cooper discretization scheme is used to ensure the positiveness and conservation of the total probability with second-order accuracy. We investigate a two-level scheme with factor-three-coarsening strategy and have a significant reduction in computations and CPU time. Numerical experiments are performed to validate the efficiency and second-order accuracy of the proposed two-level algorithm with backward time-difference schemes.
\end{abstract}
{\bf{keywords:}} Stocastic process; Fokker-Planck equation; Chang-Cooper scheme; two-level scheme; staggered grids; finite difference \\
{\bf{MSC 2010:}} 35Q84; 49K20; 65N55

\section{Introduction}
Stochastic models are used in diverse field as ecology, genetics, economics and engineering. Closed form solutions of such models are know, however, only for some of the simplest drift and diffusion functions. Therfore, there has always been a need of numerical methods that solves comples stochastic models and hence the Fokker-Planck equation. The global dynamical behaviour of a nonlinear system with noise is formally described by the probability density function (PDF) evolution along deterministic and diffusuion that satisfies the Fokker-Planck (FP) partial differential equation \cite{Risken}. The FP system that has the capability of connecting stochastic and deterministic dynamics has been applied to various applications in physics, chemistry, biology and finance \cite{Carlini, Chang, Risken, RoyBorzi2018, Spencer}.

In this article, we consider the FP equations that corresponds to the stochastic differential equations. In particular, the stochastic process defined by the following multidimensional model \cite{Mohammadi}
\begin{eqnarray}\label{stocasticmodel}
  dX_t &=& b(X_t,t)dt+\sigma(X_t,t)dW_t \\
  X(t_0) &=& X_0,
\end{eqnarray}
where $X_t\in{\mathbb{R}}^{d}$ is the state variable and $dW_t\in{\mathbb{R}}^l$ denotes the Wiener process. Moreover, $ \sigma\in{\mathbb{R}}^{d\times l}$ is a full rank dispersion matrix. Note that a statistical distribution can describes the state of the stochastic process. For this, the probability density function (PDF) distribution and the evolution of this PDF distribution can be modelled by the FP equation.

The numerical solution of the FP equations has been obtained by several researchers. One of the most popular scheme in this regard which solves the linear FP equation is the Chang-Cooper (CC) scheme introduced by Chang and Cooper in $1970$ \cite{Chang}. One of the impotent features of CC scheme is that the discrete solution preserves some intrinsic properties of the original given problem, one such as positivity and conservation of the total probability. Later, several improvements have been done \cite{Drozdov}, where we have seen high order finite difference schemes and also the nonlinear case. Finite element schemes have also been discussed, see \cite{Spencer}. It is also worth noting that some semi-analytic techniques are employed to solve the FP equation, for example, in \cite{Dehghan} the FP equation is investigated by the Adomian decomposition method. In \cite{Tatari}, variational iteration method is presented to solve the FP equation. Moreover, a finite difference scheme with cubic $C^1$-spline collocation method for solving the nonlinear Fokker-Planck equation is presented in \cite{Sepehrian}. A fast algorithm for the numerical solution of the FP equation is presented by \cite{Palleschi90, Palleschi90} and a finite difference scheme, in one-dimension, using a staggered grid to solve the Fokker-Planck equations with drift-admitting jumps is presented in \cite{Chen}. In the year 2020, the research to find the numerical solution to the stochastic models and henece the FP equation is still on; e.g., in \cite{Carlini}, a discretization scheme is developed to solve the one-dimensional nonlinear Fokker-Planck-Kolmogorov equation that preserves the nonnegativity of the solution and conserves the mass; a solution to the Fokker–Planck Equation with piecewise-constant drift is proposed in \cite{ChenChen},  a numerical method, named as information length, for measuring distances between statistical states as represented by PDF has been proposed in \cite{Anderson}. Also, there has been work on Fractional Fokker-Planck Equation as well, e.g., a space-time Petrov-Galerkin spectral method for time fractional FP equation with nonsmooth solution has been studied in \cite{Zeng} and a numerical solution of the Cauchy problem for the fractional FP equation in connection with Sinc convolution methods is proposed in \cite{Baumann&Stenger}.

In this work, we intend to solve the FP equation with linear and nonlinear drift vector and constant diffusion tensor. By doing this, the Gaussian distribution for the FP equation, which is a parabolic type differential equation that also satisfies the positiveness and conservation condition (\ref{normalizedcond_tfpe}). We present a two-level algorithm with coarsening by a factor-of-three strategy on staggered grids c.f. \cite{Butt,Butt&Borzi,Butt&Yuan} with (backward) time-difference scheme of order one and two, i.e., BDF1 and BDF2, respectively. A Chang-Cooper discretization scheme has been used to guarantee the second-order accuracy, positiveness and conservation of the total probability.

In the next Sec. {\ref{sec_tfpe}}, the Fokker-Planck equation is presented and a Chang-Cooper discretization scheme is explained in Sec. \ref{sec_disfpe}. A two-level scheme with inter-grid transfer operators is presented in Sec. \ref{sec_mg}. Numerical results are reported in Sec. \ref{sect_num} and a Sec. of conclusions is given in the last.
\section{Fokker-Planck equation \label{sec_tfpe}}
We consider the following time-dependent Fokker-Planck equation in one-dimensional computational domain $\Omega \subset \mathbb{R}$ with $Q:=\Omega \times (0,T)$ and Lipschitz boundary $\partial\Omega$:
\begin{equation}\label{tfpeq}
\frac{\partial u(x,t)}{\partial t} - \frac{\sigma^2}{2} \frac{\partial }{\partial x^2} u(x,t) + \frac{\partial}{\partial x} (f(x,t)\,u(x,t)) = 0, \qquad \mbox{ in } Q
\end{equation}
with the initial PDF distribution
\begin{equation}\label{inicond_tfpe}
u(x,0) = u_0(x), \qquad \mbox{ in }  \Omega
\end{equation}
which satisfies the positiveness and conservation of PDF distribution condition
\begin{equation}\label{normalizedcond_tfpe}
u_0 \geqslant 0, \qquad \int_{\Omega} u_0(x) \,dx = 1.
\end{equation}

The Fokker-Planck equation (\ref{tfpeq}) can be written in flux form (with non-zero source term $g$), i.e.,
\begin{equation}\label{fluxtfpe}
\frac{\partial u}{\partial t} - \nabla \cdot F(x,t) = \frac{\partial u}{\partial t} - \frac{\partial}{\partial x} F(x,t) = g(x,t)
\end{equation}
where
\begin{equation*}
\nabla  = \frac{\partial}{\partial x}, \qquad F(x,t) = B(x,t)u(x,t)+C(x,t)\frac{\partial}{\partial x} \, u(x,t)
\end{equation*}
represents the flux and the source term $g(x,t)$ has been added for the numerical investigation purposes. However, the positivity and conservation of the PDF distribution function $u(x,t)$ for the FP equation are claimed when $g(x,t)=0$. The initial condition is given by (\ref{inicond_tfpe}) and the boundary conditions are
\begin{equation}\label{bcstimefpe}
F = 0, \qquad \mbox{ on }  \partial \Omega \times (0,T).
\end{equation}

For simplicity, we choose $C(x,t)=a_{ii}(x,t)$, $B(x,t)=\partial x_i a_{ii} (x,t)$. Also, we assume that $C(x,t)$ is a positive continuous scalar function and in the case of Ornstein-Uhlenbeck process that we shall follow, $C(x,t)$ is a positive constant function and $B(x,t)$ is constant in time and linear in the spatial variable. Further, we assume that $B(x,t)$ is a function such that satisfies the Lipschitz continuity
\begin{equation*}
|B(x+h,t)-B(x,t)|\leq y \, h
\end{equation*}
where $y$ is the Lipschitz constant.
\section{Discretization on staggered grid \label{sec_disfpe}}
\begin{figure}[h!]
\begin{center}
\setlength{\unitlength}{1cm}
\begin{picture}(13,4)
\put(0.5,1.5){\framebox(1.0,1.0){$F$}}
\put(1.5,2){\line(1,0){1}}
\put(3.5,2){\line(1,0){1}}
\put(4.5,1.5){\framebox(1.0,1.0){$F$}}
\put(5.5,2){\line(1,0){1}}
\put(7.5,2){\line(1,0){1}}
\put(8.5,1.5){\framebox(1.0,1.0){$F$}}
\put(9.5,2){\line(1,0){1}}
\put(11.5,2){\line(1,0){1}}
\put(12.5,1.5){\framebox(1.0,1.0){$F$}}
\put(3,2){\circle{1}{$u$}}
\put(7,2){\circle{1}{$u$}}
\put(11,2){\circle{1}{$u$}}
\end{picture}
\caption{Staggered grid for one-dimensional FP equation} \label{staggeredgrid}
\end{center}
\end{figure}
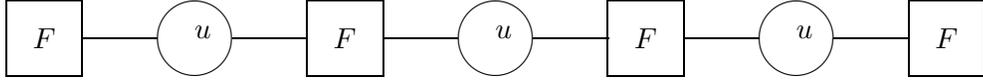

In this section, we discretize the FP equation on staggered grid, see Fig. \ref{staggeredgrid}. We use the Chang–Cooper (CC) scheme which is second-order accurate and guarantees the conservation of the total probability and positive solution to the numerical solution of FP equation, see \cite{Chang}.

We consider a one-dimensional computational domain, i.e., $\Omega = (-a,a)$. For discretization, we consider a sequence of uniform grids $\{\Omega_h\}_{h>0}$ with spatial mesh size $h$ and $N$ as the number of cells
\begin{equation*}
\Omega_h = \left\{x_i = -a+ih, i= 0,1,\ldots, N \right\} \cap \Omega
\end{equation*}
On a uniform staggered grid, the flux $F$ and PDF distribution function $u$ (solution points) are
\begin{align*}
&F_{i} = F(-a+i\,h),& &0 \leq i \leq N& \\
&u_{i} = u(-a+(i-1/2)h),& &1 \leq i \leq N.&
\end{align*}
We choose the spatial mesh size $h$ such that the boundary of the domain $\Omega$ coincide with the grid points.

Note that the time-dependent FP equation is a parabolic type equation with an additional (positiveness and conservation of PDF distribution) condition (\ref{normalizedcond_tfpe2d}). For discretization, let ${\tau}$ and $N_t$ be the time stepsize and number of time steps, respectively. We define
\begin{equation}
Q_{h,{\tau}} = \left\{ (x_i,t_m): x_i = -a+(i-1/2)h, t_m = m {\tau}; i= 1,\ldots, N; 0 \leq m \leq N_t \right\}
\end{equation}
where on $Q_{h,{\tau}}$, we mean $u_{i}^m$ the value of the grid function in $\Omega_h$ at $x_i$ and time $t_m$. The Chang-Cooper discretization scheme is used c.f. \cite{Chang} for the spatial variable and for $\frac{\partial u}{\partial t}$ approximation we use first-order backward difference formula (BDF1)
\begin{equation}\label{bdf1}
  \frac{\partial u}{\partial t} \approx \frac{u_{i}^{m+1} - u_{i}^m}{{\tau}}.
\end{equation}

The discetize version of $\nabla \cdot F$ at time $t_m$, corresponding to the time-dependent FP equation, is given by
\begin{equation}\label{tfluxh}
  \nabla \cdot F =  \frac{F_{i+1/2}^m - F_{i-1/2}^m}{h}
\end{equation}
where $F_{i+1/2}^m$ denotes the flux in $x$-direction at the point $x_i$. The discretized flux zero boundary conditions (\ref{bcstimefpe}) are given by
\begin{eqnarray}\label{bcsfpeh}
F(-1/2,t_m) = 0, \qquad F(N+1/2,t_m) = 0, \qquad \forall \, m=0,1,\cdots,N_t
\end{eqnarray}
Moreover,
\begin{eqnarray*}
F_{i+1/2}^m &=& \left[ (1-\delta_{i+1/2}^m) f_{i+1/2}^m + \frac{\sigma^2}{2h} \right]u_{i+1}^{m+1}
- \left[ \frac{\sigma^2}{2h} - \delta_{i+1/2}^m f_{i+1/2}^m \right]u_{i}^{m+1} \\
\end{eqnarray*}
and
\begin{eqnarray}\label{deltatimefpeh}
\delta_{i+1/2}^m &=& \frac{1}{\omega_{i+1/2}^m} - \frac{1}{exp(\omega_{i+1/2}^m)-1}, \qquad \omega_{i+1/2}^m = \frac{2hf_{i+1/2}^m}{\sigma^2}.\nonumber
\end{eqnarray}
Thus the discrete time-dependent FP equation $\frac{\partial u}{\partial t} - \nabla \cdot F = 0$ becomes
\begin{equation}\label{timefpeh}
\frac{u_{i}^{m+1} - u_{i}^m}{{\tau}} - \frac{F_{i+1/2}^m - F_{i-1/2}^m}{h} = 0.
\end{equation}

Conservation of the FP equation follows from the discrete FP equation and for this we use the flux form of the FP equation:
\begin{lemma}\rm
The conservation property holds
\begin{equation*}
  \sum_{i=0}^{N}\,u_{i}^{m+1}=\sum_{i=0}^{N} \,u_{i}^{m}, \qquad m\geq0.
\end{equation*}
\proof
Denote $D_t \, u_i^m = \frac{u_i^{m+1}-u_i^m}{{\tau}}$, and consider the time-dependent FP equation in flux form, i.e.,
\begin{equation*}
D_t \, u_i^m (\thickapprox \frac{\partial u}{\partial t}) = \frac{1}{h}\left(F_{i+1/2}^m-F_{i-1/2}^m\right)
\end{equation*}
in equation (\ref{bdf1}) then taking sum over $i$ gives
\begin{equation*}
\sum_{i=0}^{N} (u_i^{m+1}-u_i^m) = \frac{{\tau}}{h}\left(F_{i+1/2}^m-F_{i-1/2}^m\right).
\end{equation*}
Note that, at the boundaries, we have a zero right hand side because it is the difference of fluxes. Thus,
$$\sum_{i=0}^{N}\,u_i^{m+1}=\sum_{i=0}^{N} \,u_i^{m}.$$\endproof
\end{lemma}

For stability and convergence of the CC scheme with first-order time approximation (BDF1), we denote
\begin{eqnarray*}
  D_+ u_i &=& \frac{u_{i+1}-u_i}{h} \\
  D_- u_i &=& \frac{u_i - u_{i-1}}{h} \\
  M_\delta u_i &=& (1-\delta_{i-1})u_i + \delta_{i-1}u_{i-1}
\end{eqnarray*}
We have the following CC discretization scheme to the FP equation with non-zero source term $g$ using BDF1 for the time variable, see \cite{Chang,Mohammadi}
\begin{eqnarray*}
\frac{u_i^{m+1}-u_i^m}{{\tau}} &=& \frac{1}{h}[((1-\delta_{i}^{m})B_{i+1/2}^{m} + \frac{1}{h}C_{i+1/2}^{m})
u_{i+1}^{m+1} \\
&-& (\frac{1}{h}(C_{i+1/2}^{m}+C_{i-1/2}^{m})+(1-\delta_{i-1}^{m})B_{i-1/2}^{m} -\delta_{i}^{m}B_{i+1/2}^{m})u_{i}^{m+1} \\
&+& (\frac{1}{h}C_{i-1/2}^{m}-\delta_{i}^{m}B_{i-1/2}^{m})u_{i-1}^{m+1}]+g_j^{m+1}, \qquad i=0,1,\ldots N
\end{eqnarray*}
where
\begin{eqnarray*}
F_{i+1/2}^{m} &=& B_{i+1/2}^{m} \left((1-\delta_{i}^{m})u_{i+1}^{m+1} + \delta_{i}^{m} u_i^{m+1}\right) + C_{i+1/2}^m
\left(\frac{u_i^{m+1}-u_i^m}{h}\right)\\
\delta_{i}^m &=& \frac{1}{\omega_{i}^m} - \frac{1}{exp(\omega_{i}^m)-1} \\
\omega_{i}^m &=& \frac{h \, B_{i+1/2}^m}{C_{i+1/2}^m}
\end{eqnarray*}
with zero-flux boundary conditions, i.e., $F_{-1/2}^{m} = 0, F_{N+1/2}^{m}=0$. Note that at equilibrium the numerical fluxes must be zero, $F_{i+1/2}=0$. Therefore,
\begin{equation*}
\frac{u_{i+1}^{m+1}}{u_{i}^{m+1}} = \frac{\frac{1}{h}C_{i+1/2}^{m}-\delta_i^m B_{i+1/2}^{m}}{(1-\delta_i^m)B_{i+1/2}^{m}+\frac{1}{h}C_{i+1/2}^{m}},
\end{equation*}
and if we solve $F(x_{i+1/2},t^{m+1})=0$, we have
\begin{eqnarray*}
\frac{u_{i+1}^{m+1}}{u_{i}^{m+1}} &=& exp\left(- \int_{x_i}^{x_{i+1}}\frac{B(x,t^{m+1})}{C(x,t^{m+1})} \, dx_i\right) \\
                                  &\approx& \frac{h\,B_{i+1/2}^{m}}{C_{i+1/2}^{m}}.
\end{eqnarray*}

With this setting, the discretized FP equation with source term $g(x,t)$ is given by
\begin{equation}\label{CCFPE1d}
\frac{u_i^{m+1}-u_i^m}{{\tau}} = D_+ C_{i-1/2}^m D_-u_{i}^{m+1} + D_+B_{i-1/2}^{m}M_{\delta} u_i^{m+1} + g_i^{m+1}
\end{equation}
and for positivity, stability and convergence results, see \cite{Mohammadi}.

\section{Two-level scheme \label{sec_mg}}
In this section, we illustrate the proposed two-grid algorithm with intergrid transfer operators in details. As we know that multigrid scheme uses grids that we obtained after discretization (finite difference or finite element) and such grids are usually obtained from a coarse grid, for example, by halving the coarsest grid, see \cite{Trottenberg}. As a result, we obtained a non-nested hierarchy of grids and need extra efforts to construct intergrid transfer operators. Therefore, we note that when a coarsening by a factor-of-three is used, we obtain a nested sequence of grids. This allows us to use bilinear interpolation and straight injection and hence the implementation of intergrid transfer operators becomes easier, which we explain in details as follows.

Let $\Omega_k$ denotes the nested grids or levels with mesh size
$h_k=1/3^{k-1}$, where $k = 1,\dots,L$, and $L$ denotes the finest level. In this way, we have a variable $X_{I}^{k-1}$ at the coarse grid point $I$ on $\Omega_{k-1}$ that has the same spatial location as the variable $X_i^{k}$ at the fine grid point $i$ on $\Omega_{k}$, see Fig. \ref{staggeredgrid}
\begin{itemize}
\item  $u_{I+1/2}^{k-1}$ corresponds to $u_{i+1/2}^{k}$ for $i=3I-1$.
\end{itemize}

In the following, we present the intergrid transfer operators, that is, the prolongation and restriction operators, that are to be used in the two-level algorithm.

For prolongation, we use a quadratic interpolation using Lagrange polynomial, e.g., using a quadratic polynomial $p(x)=a_0+a_1x+a_2x^2$ that satisfies $p(x_i)=u_i, \qquad i = 1,2,3$, we have the following quadratic (second-order) Lagrange Interpolation formula
\begin{equation*}
p(x) = \frac{(x-x_2)(x-x_3)}{(x_1-x_2)(x_1-x_3)}u_1+\frac{(x-x_1)(x-x_2)}{(x_2-x_1)(x_2-x_3)}u_2+\frac{(x-x_1)(x-x_2)}{(x_3-x_1)(x_3-x_2)}u_3
\end{equation*}
Moreover, consider the space $\cU_k$ of $u^k: \Omega_k \to \R$, $k=1,\ldots,L$ such that for every two grids ${\Omega}_k$ and ${\Omega}_{k-1}$, a prolongation operator, $I^k_{k-1}:\cU_{k-1} \to \cU_{k}$ is defined which is consistent with each partition or subinterval of the discretization.
\begin{figure}[h!]
\begin{center}
\setlength{\unitlength}{1cm}
\begin{picture}(13,6)
\put(4.5,1.5){\framebox(1.0,1.0){$F$}}
\put(8.6,1.5){\framebox(1.0,1.0){$F$}}
\put(6.3,1.8){\framebox(.4,.4){$F$}}
\put(7.5,1.8){\framebox(.4,.4){$F$}}

\put(5.5,2){\line(1,0){.8}}
\put(6.7,2){\line(1,0){.8}}
\put(7.9,2){\line(1,0){.7}}

\put(5.7,2){\circle{.1}}
\put(5.9,2){\circle{.1}}
\put(6.1,2){\circle{.1}}

\put(6.9,2){\circle{.1}}
\put(7.1,2){\circle{.1}}
\put(7.3,2){\circle{.1}}

\put(8.1,2){\circle{.1}}
\put(8.3,2){\circle{.1}}
\put(8.5,2){\circle{.1}}

\put(0.5,4.5){\framebox(1.0,1.0){$F$}}
\put(1.5,5){\line(1,0){1}}
\put(3.5,5){\line(1,0){1}}
\put(4.5,4.5){\framebox(1.0,1.0){$F$}}
\put(5.5,5){\line(1,0){1}}
\put(7.5,5){\line(1,0){1}}
\put(8.5,4.5){\framebox(1.0,1.0){$F$}}
\put(9.5,5){\line(1,0){1}}
\put(11.5,5){\line(1,0){1}}
\put(12.5,4.5){\framebox(1.0,1.0){$F$}}
\put(3,5){\circle{1}{$u$}}
\put(7,5){\circle{1}{$u$}}
\put(11,5){\circle{1}{$u$}}
\end{picture}
\caption{Illustration of {\it{Straight Injection operator $I_k^{k-1}$}} on a single subintrval (partition). The the coarse grid $\Omega_{k-1}$ (upper line) after coarsening by a factor-of-three of the fine grid (bottom line) points on $\Omega_{k}$}\label{injection_stagrid}
\end{center}
\end{figure}
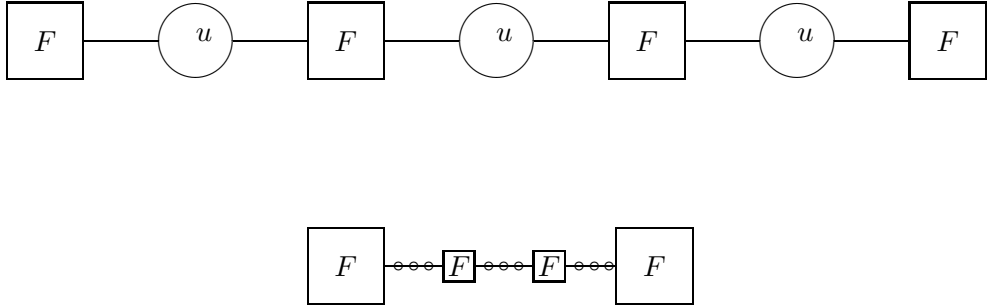

To transfer the residuals (resp. functions) from fine ${\Omega}_k$ to coarse grid ${\Omega}_{k-1}$
a straight injection, that is, $I^{k-1}_{k}: \cU_{k} \to \cU_{k-1}$ is used as a restriction operator. For example, for the flux function $F$ and PDF distribution $u$, we have the following correspondence, see Fig. \ref{injection_stagrid}
\begin{itemize}
\item  $F^{k-1}_{I} \leftarrow F^{k}_{i}$ for $i=3I-2$;
\item  $u_{I+1/2}^{k-1} \leftarrow u_{i+1/2}^{k}$ for $i=3I-1$.
\end{itemize}
Here we remark that we can use the half or full weighting as a restriction operator $I^{k-1}_{k}$. We use the straight injection because it gives a natural choice in a coarsening by a factor of three \cite{Butt&Borzi,Butt,Butt&Yuan} and hence an advantage of using the proposed multigrid scheme with factor-three coarsening.

The two-level algorithm to solve the Fokker-Planck equation is given by:
\begin{algorithm}\label{tFP_TGalgo}
{\bf TG$(m_1,m_2)$ for solving $A_k \, u_k =g_k$}.
\begin{enumerate}
\item{Set $u_k^0$ using initial condition and normalized condition (\ref{normalizedcond_tfpe});} \\

\item  Pre-smoothing:\\
$u_k^{m,(l)}=S_k(u_k^{m, (l-1)},g_k^m)$, $l=1,\dots,m_1$; \\

\item  Compute the residual $r_k^m=u_k^m - A_k (u_k^{m, (m_1)})$;
\item  Restrict the residual $r_{k-1}^m=I^{k-1}_k r_k^m$;
\item  Solve the coarse-grid (or error equation) problem on coarser level, i.e., solve
       $A_{k-1}(e_{k-1}^m)=r_{k-1}^m$ with $e_{k-1}^m :=0$ as an initial guess;
\item  Transfer the error (using interpolation operator), i.e., $I^k_{k-1}$: $e_k^m := I^k_{k-1} e_{k-1}^m$
\item  Coarse-grid correction step: $u_k^{m,(m_1 +1)}=u_k^{m,(m_1)}+e_k^m$;
\item  Apply normalized condition (\ref{normalizedcond_tfpe});
\item  {Post-smoothing on the fine grid: \\
$u_k^{m,(l)}=S_k(u_k^{m,(l-1)},g_k^m)$, $l=m_1+2,\dots,m_1+m_2+1$;}
\end{enumerate}
\end{algorithm}
\subsection{FP equation with second-order time difference \label{sec_fpe2dt2}}
In the following, we discuss the discretization of time-dependent FP equation with second-order difference scheme to the time derivative. In particular, we consider the one-dimensional time-dependent FP equation with second-order backward time difference formula (BDF2):
\begin{equation*}
\frac{\partial u(x,t)}{\partial t} \approx \frac{3u_{i}^{m+1} - 4u_{i}^{m} + u_{i}^{m-1}}{2{\tau}}
\end{equation*}
Then, we have the following discretized FP equation c.f. (Section 3.2 \cite{Mohammadi})
\begin{equation}\label{tfpeht2}
\frac{3u_{i}^{m+1} - 4u_{i}^{m} + u_{i}^{m-1}}{2{\tau}} = D_+C_{i-1/2}^{m}D_-u_i^{m+1}+D_+B_{i-1/2}^m\,M_\delta\,u_i^{m+1}+g_i^{m+1}.
\end{equation}
For conservation property, in case of one-dimensional FP equation with second-order time difference scheme (BDF2), we have by summing over $i$ and using the zero-flux boundary conditions:
\begin{equation*}
3u_{i}^{m+1} = 4u_{i}^{m} - u_{i}^{m-1}.
\end{equation*}
Then by induction and using $\sum_{i=0}^{N} \, u_i^1 = \sum_{i=0}^{N} \, u_i^0$ and $\sum_{i=0}^{N} \, u_i^m = \sum_{i=0}^{N} \, u_i^{m-1}$, we have
\begin{equation*}
\sum_{i=0}^{N} \, u_i^{m+1} = \sum_{i=0}^{N} \, u_i^m, \qquad m\geq1.
\end{equation*}
Analogously to FP equation with first-order backward time differencing BDF1 scheme, we have the conservation property for the FP equation with second-order time difference scheme (BDF2). For detailed proof about the numerical stability, convergence and positivity of the CC scheme with second-order time difference approximation scheme (BDF2), see \cite{Mohammadi}.
\section{Numerical experiments \label{sect_num}}
In this section, we present numerical examples to solve the Fokker-Plank equation with linear and nonlinear drift function to demostarte the efficiency and second-order accuract for the proposed two-level algorithm with BDF1 and BDF2, respectively. We use Matlab $2016$ on laptop $i7$, $1.86GHz$ with$4GB$ RAM, for the numerical simulations.
\subsection{Stationary FP equation}
First, we consider a stationary FP equation on $\Omega=[-6,6]$:
\begin{equation}\label{sec_fpe1d}
\frac{d}{dx} \left[-\frac{\sigma^2}{2} \frac{d}{dx}u(x) + f(x)u(x)\right] = 0.
\end{equation}
We take the diffusion coefficient $\sigma = 1$ and the linear drift function $f(x)=-x$ so that (for comparing the numerical and analyticla solution) we have an analytic solution given by $u_e = 1/exp(x^2)$. We employ the two-level Algorithm \ref{tFP_TGalgo} with $3-pre$ and $3-post$ smoothing (i.e., $m_1 = m_2 = 3$) steps. The solution error, at the discrtized level $k$, is presented in Table \ref{tab1_errfpe_exp1} based on the following discrete $L^1-norm$
\begin{equation*}
  \||u|\|_1 = h \sum_{i=1}^{N} |u_{i}|,
\end{equation*}
and discrete $L^2-norm$
\begin{equation*}
  \|u\|_2 = h^2 \sum_{i=1}^{N} u^2_{i}.
\end{equation*}
We stop the iterations when the difference of discrete $L^2-norm$ of errors of the new and old numerical approximation to $u$, i.e., when $\| u^{new} \|_2 - \| u^{old} \|_2 < tol=10^{-8}$. Number of two-grid cycles required to reach a desired tolerance with CPU time (seconds) are also reported in Table \ref{tab1_errfpe_exp1}.

\begin{figure}[ht]
\centering
\includegraphics[height=5in,width=5in]{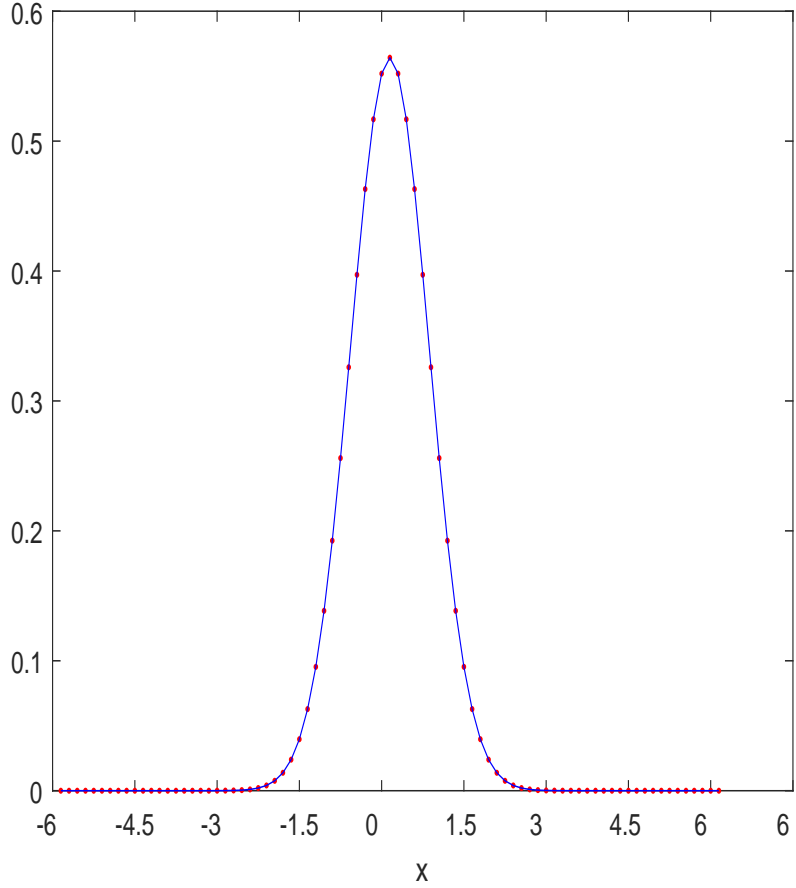}
\caption{Stationary FP equation: Numerical solution (.) and Analytical solution (-) on $N = 81$ mesh.}\label{figs_fpe1d}
\end{figure}

Further, the numerical and analytical solution for this FP model is depicted in Fig. \ref{figs_fpe1d} on $N=81$. This shows that the proposed algoritm track the desired PDF, that is, for the stationary case, we have the match of numerical and analytical PDF distributions.
\begin{table}[h!]
\caption{Error history for stationary FP equation}
\begin{center}
\begin{tabular}{lcccc}
\hline
$N$ & $\||u - u_e\||_1$ & $\|u - u_e\|_2$ & \#TG & CPU\\
\hline
$27$ & $1.7766e-09 $ & $1.0398e-09$ & $10$ & $0.07$\\
$81$  & $3.3659e-10$  & $1.6667e-10$ & $09$ & $0.09$\\
\hline
\end{tabular}
\end{center}\label{tab1_errfpe_exp1}
\end{table}

\subsection{Time-dependant FP equation with linear drift}
In this section, we show the second-order convergence of the proposed two-level scheme with BDF1 (resp. BDF2). For this, we consider a time-dependent Fokker-Planck equation (\ref{tfpeq}) given by \cite{Roy&Borzi2017} with initial and boundary conditions on $Q:=\Omega \times (0,T)$ with $\Omega=[-6,6]$ and $T =1$, i.e., we consider the FP equation (\ref{tfpeq}) within the framework of Ornstein-Uhlenbeck process \cite{Chen,Risken} with analytic solution:

The source term $g(x,t) = 1/exp(x^2+t)$ and drift
\begin{equation*}
f_1(x,t) = -x
\end{equation*}
and diffusion coefficient $\sigma = 1$ which results in exact solution
\begin{equation*}
  u_e(x,t)= 1/exp(x^2+t).
\end{equation*}
We employ the two-grid Algorithm \ref{tFP_TGalgo} with $3-pre$ and $3-post$ smoothing steps. Results for the discrete $L^1-norm$ and $L^2-norm$ of errors are reported in Table \ref{tab_timesourcefpe1d}), where $$\||u|\|_1 = h^2 \, {\tau}\sum_{m=0}^{N_t} \sum_{i=1}^{N} |u^k_{i}|$$
which we identify with $L^2_{{\tau}}(0,T;L^1)$ and discrete $L^2$-norm $\|\cdot\|_{L^2_{h,{\tau}}(Q)}$ , i.e.,
\begin{equation*}
  \|u\|_2 = {\tau} \, h^2 \sum_{i=1}^{N} |u^2_{i}|.
\end{equation*}

We take the time step size ${\tau} = (\frac{1}{81})(\frac{T}{3^L})$, where $L$ denote the fine level. Second-order accuracy $O(h^2+{\tau})$ is observed in the numerical results for the proposed two-level scheme, see Table \ref{tab_timesourcefpe1d}. In fact, we have a reduction in errors by a factor of nine (i.e., a factor $3^2$) as we refine the mesh by factor of $3$ for the $L^2-norm$ of errors.

In Table \ref{tab_time2sourcefpe1d}, we report results of discrete $L^1-norm$ and $L^2-norm$ of errors with CUP time (second) with second-order time differencing scheme (BDF2) using the proposed two-level scheme to the same numerical example. We use the two-level scheme at $t=t_1$ with first-order time difference scheme $BDF1$. Second-order accuracy $O(h^2+{\tau})$ is achieved , that is, we have a reduction in errors by a factor of $3^2$ as we refine the mesh by factor of $3$ for the $L^2-norm$ of errors.
\begin{table}[h!]
\caption{Convergence of FP equation with BDF1}
\begin{center}
\begin{tabular}{lcccc}
\hline
$N \times N_t$ & $\||u - u_e\||_1$ & $\|u - u_e\|_2$ & $\#TG$ & $CPU$\\
\hline
$81 \times 81$ & $6.2447e-05$ & $1.9392e-06$ & 4 & 0.14 \\
$243 \times 243$& $2.3321e-05$ & $2.4187e-07$ & 2 & 0.22 \\
$729 \times 729$ & $4.1551e-06$ & $1.6076e-08$ & 2 & 0.84 \\
$2187 \times 2187$ & $1.2060e-06$ & $1.4322e-09$ & 2 & 6.13 \\
\hline
\end{tabular}
\end{center}\label{tab_timesourcefpe1d}
\end{table}

\begin{table}[h!]
\caption{Convergence of FP equation with BDF2}
\begin{center}
\begin{tabular}{lcccc}
\hline
$N \times N_t$ & $\||u - u_e\||_1$ & $\|u - u_e\|_2$ & $\#TG$ & $CPU$\\
\hline
$81 \times 81$ & $6.2450e-05$ & $1.9393e-06$ & 4 & 0.08 \\
$243 \times 243$& $2.3330e-05$ & $2.4195e-07$ & 2 & 0.33 \\
$729 \times 729$ & $4.1450e-06$ & $1.6049e-08$ & 2 & 0.89 \\
$2187 \times 2187$ & $1.2077e-06$ & $1.4330e-09$ & 2 & 6.16 \\
\hline
\end{tabular}
\end{center}\label{tab_time2sourcefpe1d}
\end{table}

\begin{figure}[ht]
\centering
\includegraphics[height=2.5in,width=2.5in]{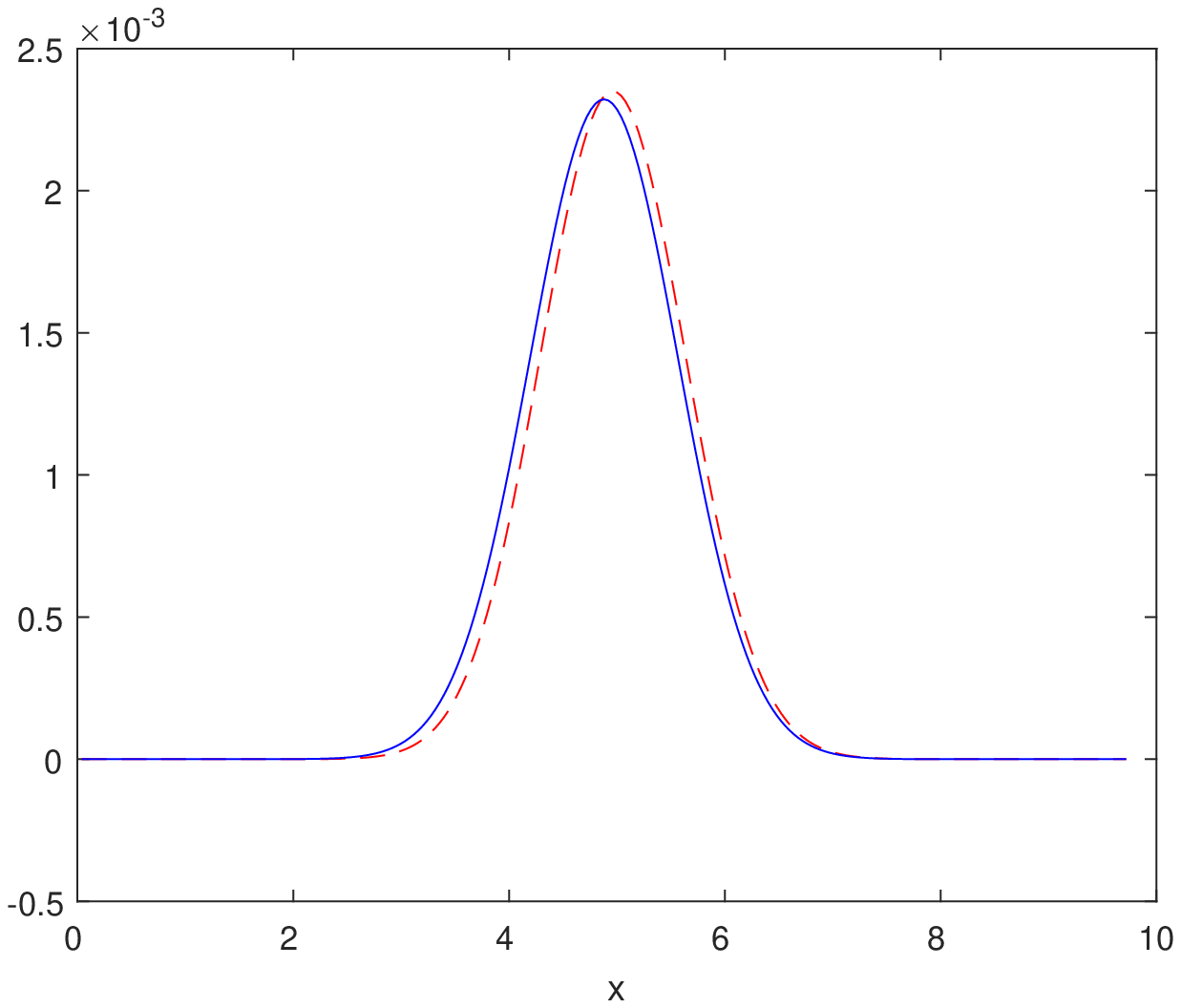}
\includegraphics[height=2.5in,width=2.5in]{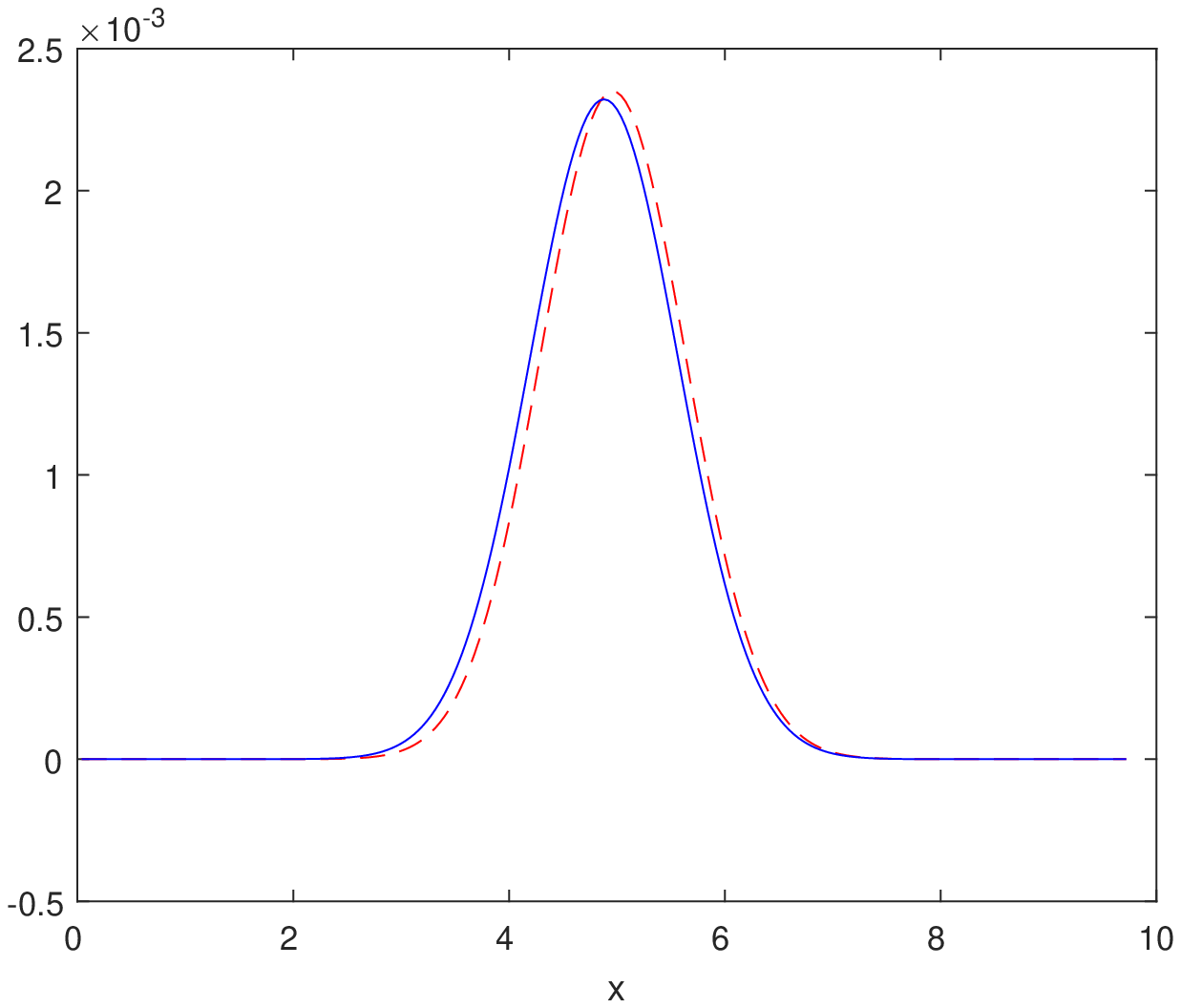}
\caption{Time-dependent FP equation \cite{Mohammadi}: Numerical solution (dashes) and Analytical solution (solid line) on $N = 243$ with BDF1 (left) and BDF2 (right), respectively.}\label{figs_fpe1dMoh}
\end{figure}

Next, to have a comparison of the proposed two-level scheme with the Chang-Cooper with first-order time backward difference (CC-BDF1) scheme given by Mohammadi and Borzi \cite{Mohammadi}, we consider the following (Ornstein-Uhlenbeck process) FP equation in $Q:=\Omega \times [0,T]$ and take $B(x,t)=x$, $C(x,t) = \sigma^2$:
\begin{equation}
\frac{\partial u(x,t)}{\partial t} = \partial_x \, \left(B(x,t) \, u(x,t) + C(x,t)\, \partial_x \, u(x,t)\right) + g(x,t), \qquad \mbox{ in } Q
\end{equation}
where the source term is given by
\begin{equation*}
g(x,t)=\frac{(a-x)(2x-a)}{exp((x-a/2)^2+t)}
\end{equation*}
and the drift function $f(x,t) = x$. In partcular, we choose $\Omega = [0,a]$ with $a= 10$ and $\sigma = 1$, $T=1$. Furthermore, the initial condition is given by
\begin{equation*}
u(x,0) = \frac{1}{exp((x-a/2)^2)}
\end{equation*}
with flux zero boundary and the exact solution is given by
\begin{equation*}
u_e(x,t)= \frac{1}{exp((x-a/2)^2+t)}.
\end{equation*}
We employ the two-grid Algorithm \ref{tFP_TGalgo} with $3-pre$ and $3-post$ smoothing steps. Results for the discrete $L^2-norm$ of errors are reported in Table \ref{tab_tfpe1d_expMohBor}. We take the same time step size ${\tau} = 0.01(\frac{1}{3^L})^2$, where $L$ is the fine level in the two-grid Algorithm. Second-order accuracy $O(h^2+{\tau})$ is observed in the numerical results for the proposed two-level scheme, see Table \ref{tab_tfpe1d_expMohBor}. Moreover, we present the numerical results of Chang-Cooper with first-order time difference (BDF1) scheme given by Mohammadi and Borzi \cite{Mohammadi}, in Table \ref{tab_tfpe1d_ComExpMohBor}. From Table \ref{tab_tfpe1d_expMohBor}-\ref{tab_tfpe1d_ComExpMohBor}, we clearly see that our proposed scheme gives better accuracy as compared to the numerical results given by \cite{Mohammadi}. Relative discrete $L^2_h$-norm of errors, on $N=81,N=243,N=729$ are recorded as $1.2253e-1, 8.2669e-2$ and $6.3865e-2$, respectively. In Fig \ref{figs_fpe1dMoh}, we depict the numerical and analytic solution on $N=243$ at $T=1$ with BDF1 and BDF2 to showcase the accuracy of the proposed Two-Level Algorithm.
\begin{table}[h!]
\caption{Convergence of proposed two-level with BDF1 for numerical example c.f. \cite{Mohammadi}}
\begin{center}
\begin{tabular}{lcccc}
\hline
$N \times N_t$  & $\|u - u_e\|_2$ & $\#TG$ & $CPU$\\
\hline
$81 \times 81$ &  $4.6050e-09$ & 2 & 0.16 \\
$243 \times 243$& $1.1507e-10$ & 2 & 0.27 \\
$729 \times 729$& $3.2926e-12$ & 2 & 0.96 \\
\hline
\end{tabular}
\end{center}\label{tab_tfpe1d_expMohBor}
\end{table}

\begin{table}[h!]
\caption{Convergence of proposed two-level with BDF2 for numerical example c.f. \cite{Mohammadi}}
\begin{center}
\begin{tabular}{lcccc}
\hline
$N \times N_t$  & $\|u - u_e\|_2$ & $\#TG$ & $CPU$\\
\hline
$81 \times 81$ &  $4.6051e-09$ & 2 & 0.11 \\
$243 \times 243$& $1.1508e-10$ & 2 & 0.31 \\
$729 \times 729$& $3.2930e-12$ & 2 & 0.91 \\
\hline
\end{tabular}
\end{center}\label{tab_t2fpe1d_expMohBor}
\end{table}

\begin{table}[h!]
\caption{CC-BDF1 scheme for FP equation by Mohammadi and Borzi \cite{Mohammadi}}
\begin{center}
\begin{tabular}{lccc}
\hline
$N \times N_t$   & $\|u - u_e\|_2$\\
\hline
$50 \times 50$   & $1.34e-2$\\
$100 \times 100$ & $3.50e-3$\\
$200 \times 800$ & $8.80e-4$\\
\hline
\end{tabular}
\end{center}\label{tab_tfpe1d_ComExpMohBor}
\end{table}

\subsection{FP equation with nonlinear drift \label{exp_tnonlinearfpe}}
\begin{figure}[ht]
\centering
\includegraphics[height=2.5in,width=2.5in]{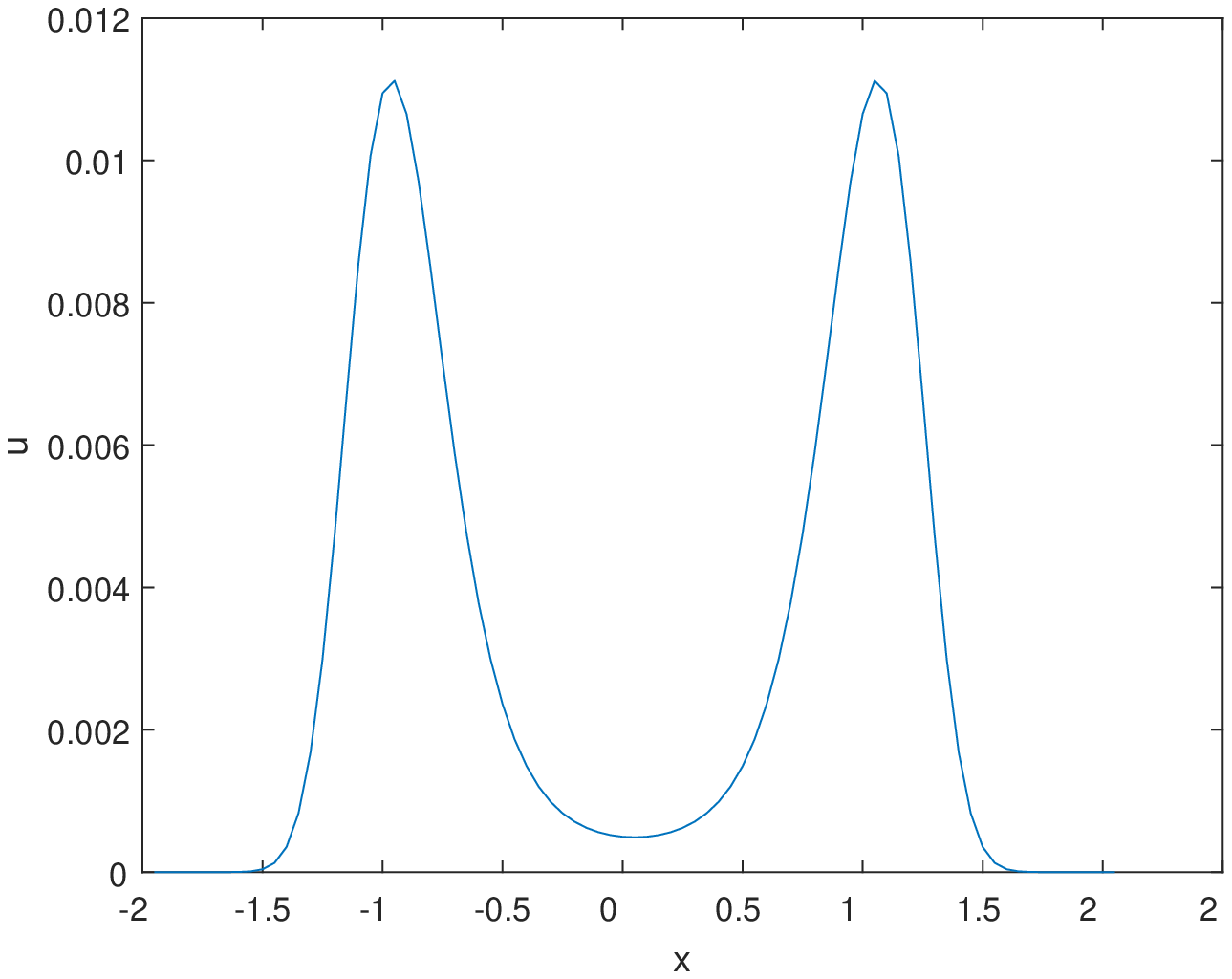}
\includegraphics[height=2.5in,width=2.5in]{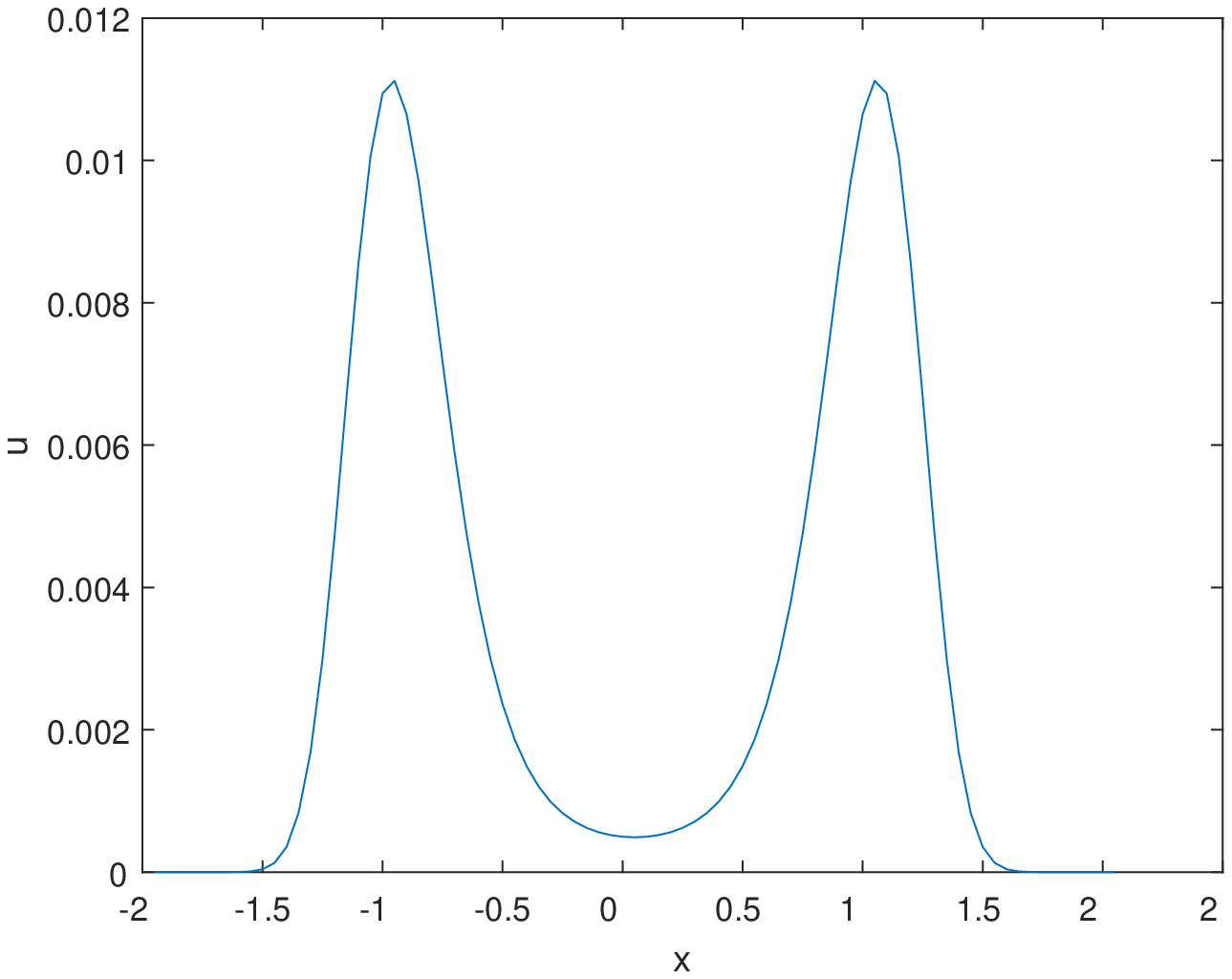} \\
\includegraphics[height=2.5in,width=2.5in]{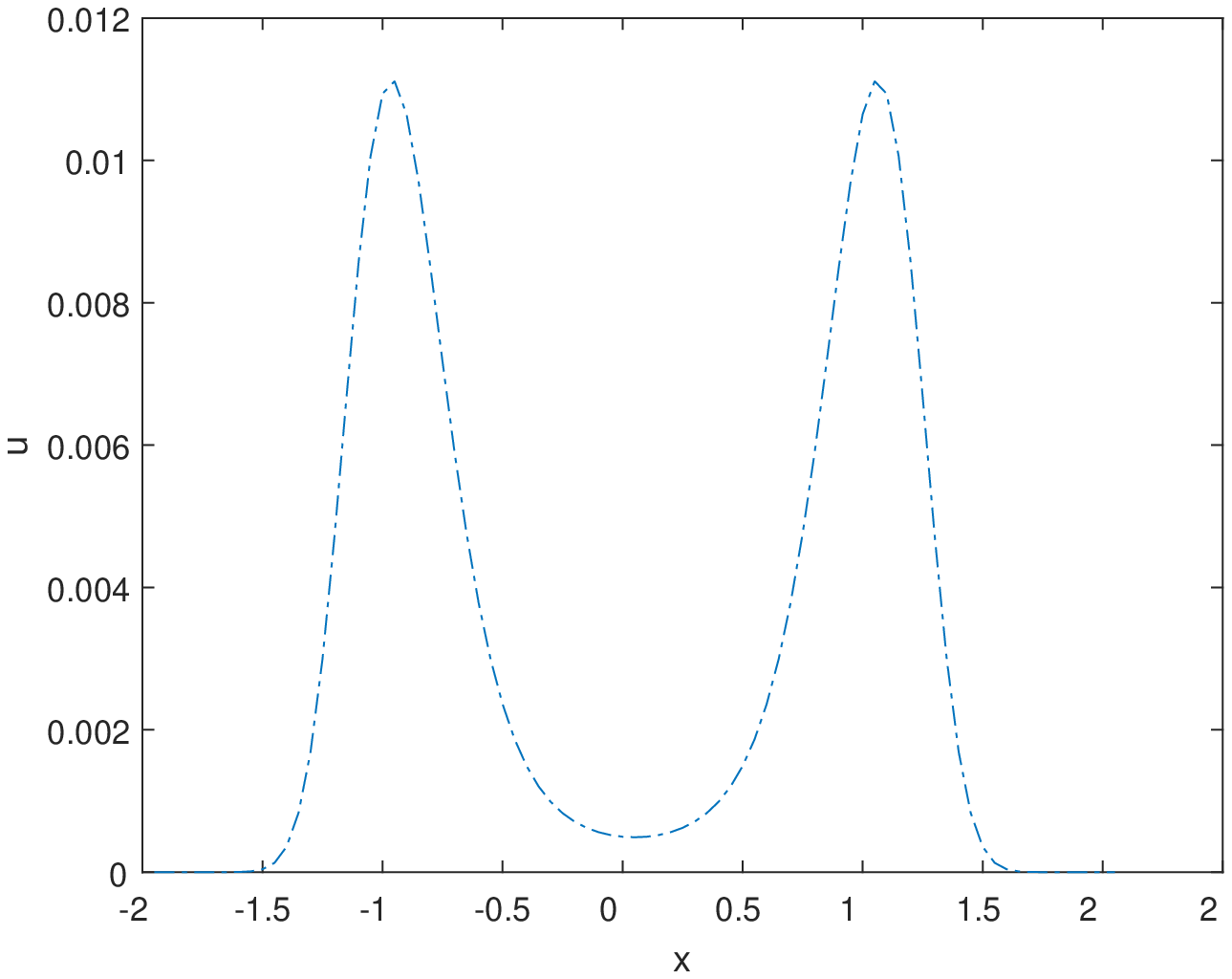}
\includegraphics[height=2.5in,width=2.5in]{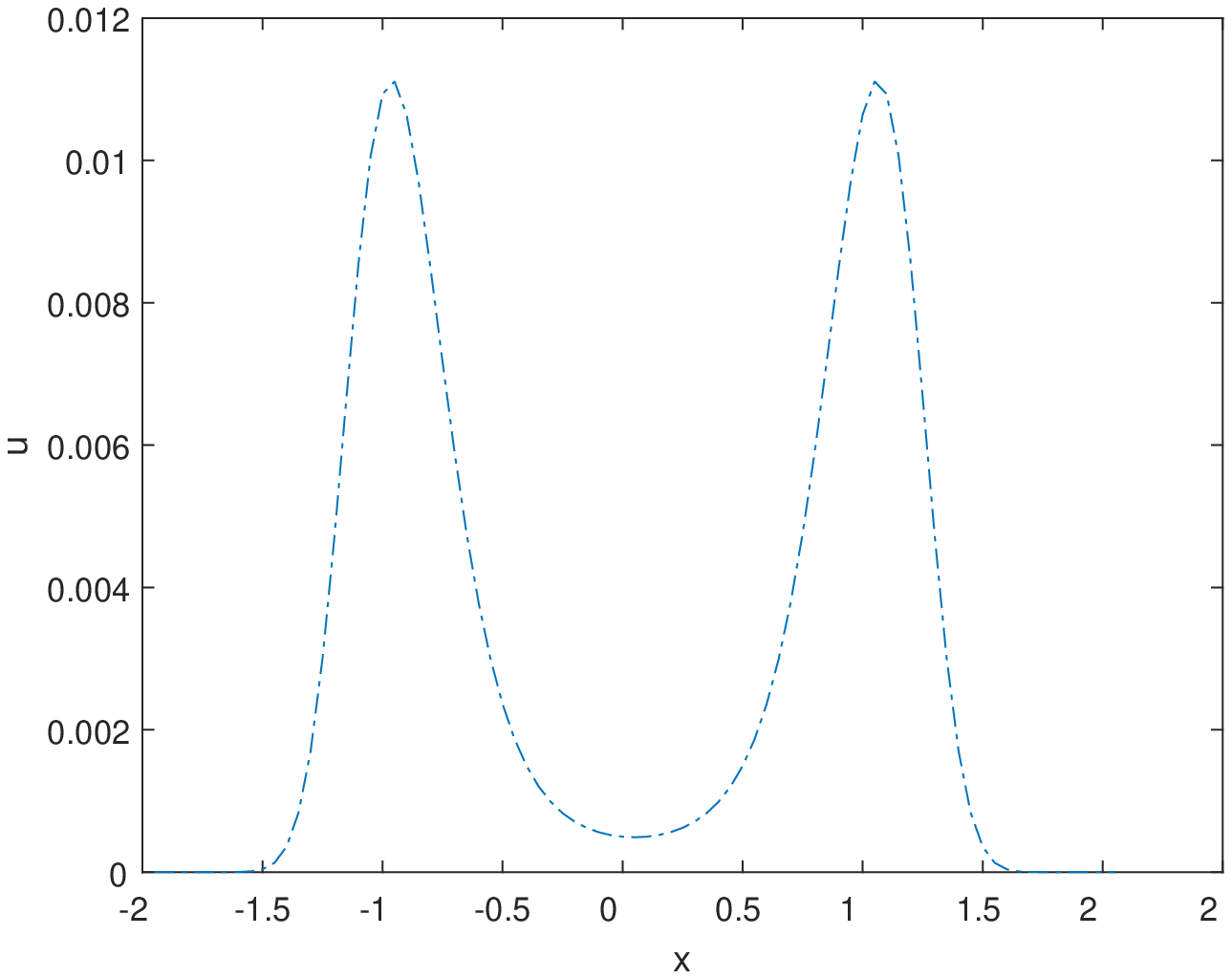} \\
\includegraphics[height=2.5in,width=2.5in]{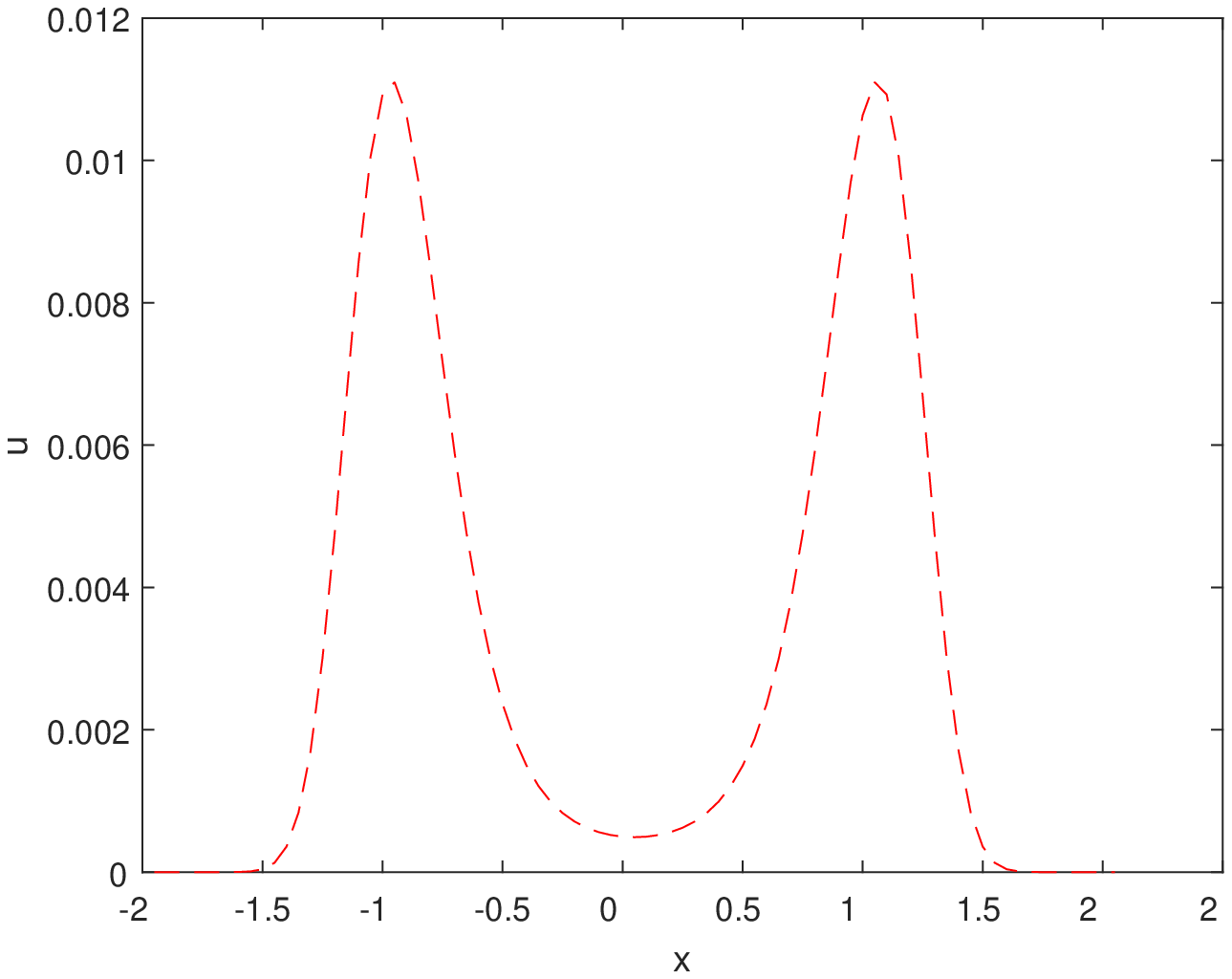}
\includegraphics[height=2.5in,width=2.5in]{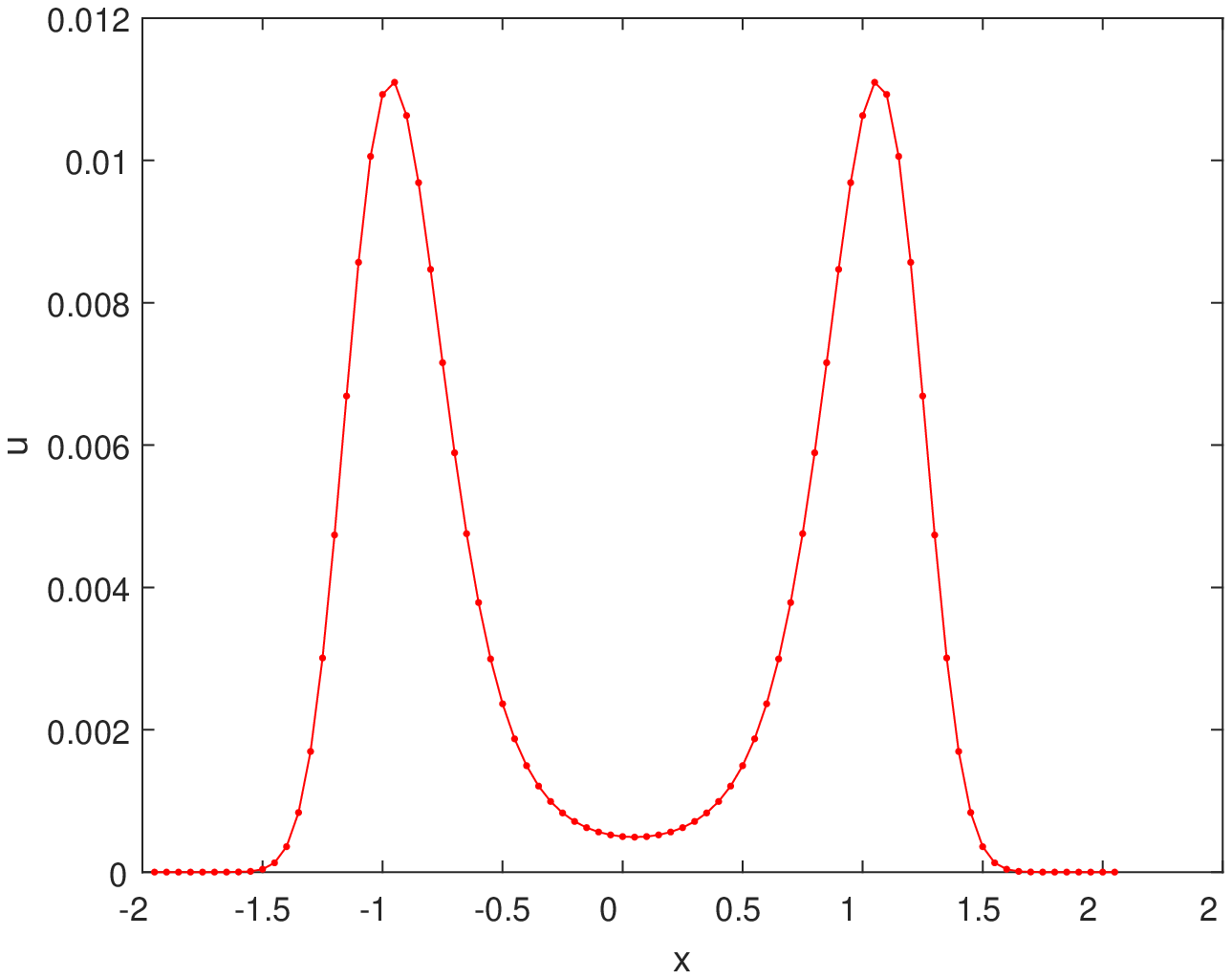} \\
\caption{Nonlinear FP equation: Numerical solution when T = 0.5 (top left); T = 1 (top right); T = 3 (middle left); T = 5 (middle right);T = 15 (bottom left); T = 30 (bottom right); on $N = 81 = N_t$ mesh, i.e., with $h=1/81$, and  ${\tau}=(1/81)(T/81)$.}\label{figs_nonlineartfpe1d_N81}
\end{figure}

In the following, we consider the nonlinear process given by Harrison \cite{Harrison}
\begin{equation}\label{nonlinearStochastic}
dX = (X-X^3)\,dt + \sigma \, dW
\end{equation}
with the corresponding FP equation
\begin{equation}\label{tnonlfpeq}
\frac{\partial u(x,t)}{\partial t} - \frac{\sigma^2}{2} \partial^2_{xx} u(x,t) + \partial_{x} (f(x,t)\,u(x,t)) = 0, \qquad \mbox{ in } Q \end{equation}
with the initial PDF distribution (\ref{inicond_tfpe}) and drift
\begin{equation*}
f(x) = x-x^3.
\end{equation*}
Moreover, we take the diffusion coefficient $\sigma = 0.4$. An analytic solution to (\ref{tnonlfpeq}) is not known. However, the steady state solution is given by
\begin{equation}\label{anasol}
u(x) = C exp((x^2-0.5x^4)/\sigma^2)
\end{equation}
where $C$ is the normalized constant. The numerical solution to this FP equation with the nonlinear drift function is depicted for $T = 0.5, 1.0, 3.0, 5.0, 15.0$ and $T=30$, respectively, on $N=81$ mesh. Moreover, to have a comparison with the results given by \cite{Harrison}, where a numerical solution of the FP equation using {\it{moving finite elements}} is presented, see Fig. \ref{figs_nonlineartfpe1d_N81}. As given by \cite{Harrison}, the deterministic equation $dx/dt=x-x^3$ has two asymptotically stable equilibria at $x = 1$ and $x = -1$ which can been seen in Fig. \ref{figs_nonlineartfpe1d_N81}. In addition, we have a symmetric bimodal distribution as a result of our numerical two-level scheme which is presented in \cite{Harrison}.
\section{Conclusions \label{sect_con}}
A two-level scheme with coarsening by a factor-of-three strategy was proposed to solve the Fokker-Planck equation with linear (nonlinear) drift function. The Chang-Cooper scheme was used to discretize the FP equation on staggered grids. Second-order accuracy, that is, $O(h^2+\tau)$ and $O(h^2+\tau^2)$ , was achieved in the numerical results using second-order differences for the spatial variable and first-order (resp. second-order) time differences BDF1 (resp. BDF2). Results of numerical examples outperform the existing numerical works on FP equation (in particular to the Ornstein-Uhlenbeck process given by \cite{Mohammadi} and to the nonlinear FP equation \cite{Harrison}). A natural extension to two-dimensional FP equation is under investigation and more complicated FP equations, that is, FP equations with nonlinear drift functions that depend on space and time both, is also our future work.

\end{document}